\documentclass[11pt]{article}
\RequirePackage{amsthm,amsmath,graphicx,epsfig,amssymb,comment}
\RequirePackage[numbers]{natbib}
\RequirePackage[colorlinks,citecolor=blue,urlcolor=blue]{hyperref}
\usepackage[utf8]{inputenc}
\usepackage[T1]{fontenc}
\usepackage{lmodern}
\newtheorem{theorem}{Theorem}
\newtheorem*{definition}{Definition}

\newtheorem*{rem}{Remark}
\newtheorem{proposition}{Proposition}
\newtheorem{cor}{Corollary}
\title{Multi-level Bayes and MAP monotonicity testing}
\author{Golubev Yu.  \thanks{ Institute for Information Transmission Problems, Moscow, Russia.  
\texttt{e-mail: golubev.yuri@gmail.com}}
\ and Pouet C.
\thanks{Centrale Marseille, Aix Marseille Univ, CNRS, 
 I2M,  Marseille,  France.
  \texttt{e-mail: cpouet@centrale-marseille.fr}}}
\begin{document}
\maketitle
\begin{abstract}
In this paper, we  develop Bayes and maximum  a posteriori probability (MAP)  approaches to monotonicity testing. In order to simplify this problem, we consider  
a simple white Gaussian noise model and with the help of the Haar transform we  reduce it to the equivalent problem of testing   positivity  of the Haar  coefficients. 
This approach permits, in particular, to understand  links between monotonicity testing  and sparse vectors detection,  to construct   new tests, and to prove their optimality without supplementary assumptions. 
The main idea in  our construction of  multi-level tests is based on some invariance properties  of specific 
probability distributions.  Along with Bayes and MAP tests, we construct  also  adaptive
multi-level   tests that are free from the prior information about the sizes of non-monotonicity segments of the function.
\end{abstract} 


\section{Introduction}

The literature on non-parametric  monotonicity testing deals usually  with the model
\[
Y=f(X)+\xi,
\]
where $Y$ is a scalar dependent random variable, $X$ a scalar independent random variable, $f(\cdot)$ 
an unknown function, and $\xi$ an unobserved scalar random variable with $\mathbf{E}\{\xi|X\} = 0$.
  We are interested in testing the null hypothesis, $\mathbf{H}_0$ that $f(x)$ is increasing against the alternative, 
  $\mathbf{H}_1$ that there are $x_1$ and $x_2$ such that $x_1 < x_2$ and $f(x_1) > f(x_2)$. 
The decision is to be made based on the i.i.d. sample $\{X_i, Y_i\}_{1\le i\le n}$ from the distribution of $(X, Y )$. Typical applications of monotonicity testing are related to econometric
 models, see, e.g., Chetverikov \cite{C}.

Usual approaches to this problem have  in their core  simple heuristic ideas and assumptions.  So,  the tests proposed in Gijbels et. al. \cite{GH} and Ghosal, Sen, and van der Vaart \cite{GS} 
are based on the signs of $(Y_{i+k} -Y_{i})(X_{i+k} -X_{i})$. Hall and Heckman \cite{HH} developed a test based 
on the slopes of local linear estimates of $f(\cdot)$. Along with these papers we can cite 
 Schlee \cite{S}, Bowman, Jones, and Gijbels 
\cite{BJ}, D\"{u}mbgen and Spokoiny \cite{DS}, Durot \cite{Du}, Baraud, Huet, and Laurent \cite{BH},  
Wang and Meyer \cite{WM},
and Chetverikov \cite{C}. As to typical hypothesis  about $f(\cdot)$, it is often assumed  that $f(x)$ is a Lipschitz  function, i.e.,
\[
|f(y)-f(x)|\le L |y-x|,
\]
where the constant $L<\infty$ may be known or unknown.

In this paper,  we look at the problem of monotonicity testing from a little different and less intuitive viewpoint. As we will see below, our approach permits, in particular, to understand  links between this problem and sparse vectors detection and  to construct   new powerful tests. 
  In order to simplify technical details and to get rid of supplementary assumptions,  we begin with   monotonicity testing of an unknown function $f(t),\, t \in [0,1]$, in the so-called white noise model similar to that one considered in 
\cite{DS}. So, it is assumed we have at our disposal the  noisy data
\begin{equation}\label{eq1}
Y(t)=f(t)+\sigma n(t), \ t \in [0,1],
\end{equation} 
where $n(\cdot)$ is a standard white Gaussian noise and $\sigma>0$ is a known noise level. With  the help of  these observations
we  want to test
\begin{align*}
&\textsl{the null hypothesis}\\
&{\bf H}_0:\  f'(t)\ge 0,\,  \text{for all}\ t\in [0,1], \\
&\textsl{vs.   the alternative}\\
&{\bf H}_1:\ f'(t)< 0,\,  \text{for some}\ t\in [0,1].
\end{align*}

Our approach  to this  problem  is based on  estimating the following linear functionals: 
\begin{equation*}
\theta_{h,t}(f)\stackrel{\rm def}{=}\frac{1}{{h}}\int_t^{t+h} f(u)\, du- \frac{1}{{h}} \int_{t-h}^t f(u)\, du
\end{equation*}
for  all $h,t$ that   are admissible, i.e., such that $[t-h,t+h]\subseteq [0,1]$. 
It is clear that   $\theta_{h,t}(f)/h$ may be interpreted as   approximations of  the derivative $f'(t)$ since 
\[
\lim_{h\rightarrow 0} \frac{\theta_{h,t}(f)}{h}=f'(t),
\]
for any given $t\in (0,1)$.

With the help of   \eqref{eq1}, the functionals  $\theta_{h,t}(f)$ are  estimated   as follows:
\[
\hat{\theta}_{h,t}(Y) =\frac{1}{{h}}\int_t^{t+h} Y(u)\, du-  \frac{1}{{h}} \int_{t-h}^t Y(u)\, du 
\]
and these estimates admit  the obvious representation
\begin{equation}\label{eq2}
\hat{\theta}_{h,t}(Y)={\theta}_{h,t}(f)+\sigma_h\xi_{h,t},
\end{equation}
where 
\[
\sigma_{h}=\sigma\sqrt{\frac{2}{h}}, \quad \xi_{h,t} =\frac{1}{{\sqrt{2h}}}\biggl[\int_t^{t+h} n(u)\, du- \int_{t-h}^t n(u)\, du\biggr]\sim \mathcal{N}(0,1).
\]

Notice  that if $\mathbf{H}_0$ is true, then  
$\theta_{h,t}(f)\ge 0  $ for all admissible $ h,t$, otherwise
 ($\mathbf{H}_1$ is true)  there exist $h',t'$ such that
$
\theta_{h',t'}(f)<0$. That is why in what follows we will focus on testing
\begin{equation}\label{eq3}
\begin{split}
&\textsl{the null hypothesis}\\
&{\bf H}_{0}:\ \theta_{h,t}(f) \ge 0,\ \textit{for all  admissible}\ h,t\\
&\textsl{vs.   the alternative}\\
&{\bf H}_{1}:\ \theta_{h,t}(f)< 0, \  \textit{for some admissible}\ h,t
\end{split}
\end{equation}
  based on the observations  \eqref{eq2}.

 Let us denote for brevity
\[
\theta_{h,t}=\theta_{h,t}(f), \quad  \hat{\theta}_{h,t}=\hat{\theta}_{h,t}(Y) .
\] 

In order to explain our approach to the problem \eqref{eq3}, we begin with  the simple case assuming that $h,t$ are given.  So,  we have to  test  two
 composite hypotheses $$ \mathbf{H}_{0}^{h,t}: \theta_{h,t} \ge 0 \ \text{ vs.} \ \mathbf{H}_{1}^{h,t} :\theta_{h,t} < 0.$$

Intuitively,  the most powerful   test with the type I error probability $\alpha$  rejects ${\bf H}_{0}^{h,t}$  if
\begin{equation}\label{eq4}
\hat{\theta}_{h,t}\le -\sigma_h t_\alpha,
\end{equation}
where $t_\alpha$ is $\alpha$-value  of the standard Gaussian distribution, i.e., a solution to
\[
\Phi(t_\alpha)=1-\alpha,
\]
where
\[
\Phi(x)=\frac{1}{\sqrt{2\pi}}\int_{-\infty}^x \exp\biggl(-\frac{x^2}{2}\biggr)\, dx.
\]

Of course, there exist  a lot of motivations for this test. In this paper, we 
 make use of the so-called improper Bayes approach 
 assuming that $\theta_{h, t}$ in \eqref{eq2} is  a random variable uniformly distributed on the interval $[0,A], \, A>0$,   if $\mathbf{H}_{0}^{h,t}$  is true, and on $[-A,0]$  if $\mathbf{H}_{1}^{h,t}$  is true.
 So,  we observe a random variable $\hat{\theta}_{h,t} $ with the probability density 
 \[
p_0^A(x|  \mathbf{H}_{0}^{h,t}\ \text{is true})=  \frac{1}{A}\int_{0}^A \exp\biggl[-\frac{(x-\theta)^2}{2\sigma_h^2}\biggr]d\theta 
 \]
 and 
 \[
p_1^A(x| \mathbf{H}_{1}^{h,t}\ \text{is true})=  \frac{1}{A}\int_{-A}^0 \exp\biggl[-\frac{(x-\theta)^2}{2\sigma_h^2}\biggr]d\theta. 
 \] 
Thus, we deal with  the simple hypothesis testing and by the Neyman-Pearson lemma, the most powerful test at significance level $\alpha$  rejects $\mathbf{H}_{0}^{h,t}$
when
\[
\frac{p_1^A(\hat{\theta}_{h,t})}{p_0^A(\hat{\theta}_{h,t})}\ge t^A_\alpha.
\]
Taking the limit in this equation as $A\rightarrow\infty$, we arrive at the improper  Bayes test  that 
rejects $\mathbf{H}_{0}^{h,t}$ if 
\begin{equation}\label{eq5}
\begin{split}
S\biggl(\frac{\hat{\theta}_{h,t}}{\sigma_h}\biggr)
\ge t'_\alpha,
\end{split}
\end{equation}
where  
\begin{equation}\label{eq6}
\begin{split}
\displaystyle
S(x)=&
\frac{\displaystyle\int_{-\infty}^ {0}  \exp\bigl[-{(x-\theta)^2}/{2}\bigr] \, d\theta}{\displaystyle
\int_{0}^\infty\exp\bigl[-{(x-\theta)^2}/{2}\bigr]\, d\theta}
=\frac{1}{\Phi(x)}-1.
\end{split}
\end{equation}
Since $S(x)$ is   decreasing  in $x\in \mathbb{R}$, the tests \eqref{eq4} and \eqref{eq5}  are obviously equivalent. 

In what follows,  we will make use of the following asymptotic result:
\begin{align}\label{eq8}
S(x)=&\biggl[1+O\biggl(\frac{1}{x^{2}}\biggr)\biggr]\sqrt{2\pi}(1-x)\exp\biggl(\frac{x^2}{2}\biggr), \ \text{as}\ x\rightarrow-\infty.%\\
%S(x)=&\biggl[1+O\biggl(\frac{1}{x^{2}}\biggr)\biggr]\frac{1}{\sqrt{2\pi}(x+1)}\exp\biggl(-\frac{x^2}{2}\biggr), \ \text{as}\ x\rightarrow \infty.\nonumber
\end{align}

Along with this method, one can apply the maximum likelihood (ML) or   minimax approaches.
Finally, all these methods result  in \eqref{eq4} but their initial forms  are different. For instance, the  ML test rejects  $\mathbf{H}_{0}^{h,t}$ when
\begin{equation}\label{eq7}
\begin{split}
\frac{ \displaystyle\max_{\theta<0} \exp\bigl\{-(\hat{\theta}_{h,t}-\theta)^2/(2\sigma_h^2)\bigr\} }{ 
\displaystyle
  \max_{\theta>0} \exp\bigl\{-(\hat{\theta}_{h,t}-\theta)^2/(2\sigma_h^2)\bigr\}}
  = \exp\biggl\{-\frac{\hat{\theta}^2_{h,t}}{2\sigma_h^2}{\rm sign}(\hat{\theta}_{h,t})\biggr\} \ge t''_\alpha.
  \end{split}
 \end{equation}
 
Emphasize that from a viewpoint of  testing  $\mathbf{H}_{0}^{h,t}$ vs. $\mathbf{H}_{1}^{h,t}$ 
there is no difference between   \eqref{eq7} and \eqref{eq5}, but the aggregation of  these methods  for testing $\mathbf{H}_{0}$ vs. $\mathbf{H}_{1}$ from \eqref{eq3} results in   different tests. In this paper, we make use of  the tests defined by
\eqref{eq5} since their aggregation is  simple.

In order to aggregate the statistical tests, we will make use of the so-called multi-resolution approach assuming  that
\begin{enumerate}
\item \textit{$h$ belongs to the following set of dyadic bandwidths}
\[
\mathcal{H}\stackrel{\rm def}{=}\biggl\{\frac{1}{2},\frac{1}{4}, \ldots \frac{1}{2^{k}},\ldots \biggr\};
\] 
\item \textit{$t$ belongs to the family of dyadic grids $\mathcal{G}_h,\, h\in \mathcal{H}$, defined by}
\[
\mathcal{G}_h\stackrel{\rm def}{=}\bigl\{ h, 3h,\ldots, 1-h\bigr\}, \ h\in \mathcal{H}.
\]
\end{enumerate}
There are  simple arguments motivating  these assumptions
\begin{itemize}
\item   random variables $\xi_{h,t}$ and $\xi_{h',t'}$ in \eqref{eq2} are independent if  $\{h,t\}\neq \{h',t'\}$. This  fact
simplifies significantly the statistical analysis of  tests. 
\item $\sqrt{h/2} \, \hat{\theta}_{h,t}$ are the Haar coefficients admitting a  fast   computation in the discrete version of \eqref{eq1}.
\end{itemize}

%A Bayes multi-resolution test has the following  simple structure. 

\section{Testing at a given resolution level} 
Let us fix some  bandwidth $h\in \mathcal{H}$ and denote 
for brevity by $n_h=1/(2h)$.  In this section,  we focus on  testing %(based on the observations from \eqref{eq2})
\begin{align*}
&\textsl{the null hypothesis}\\
&{\bf H}_{0}^h:\ \theta_{h,t} \ge 0\ \text{for all}\ t\in \mathcal{G}_h\\
&\textsl{vs. the alternative}\\
&{\bf H}_{1}^h:\ \theta_{h,t}< 0 \ \text{for some}\ t\in \mathcal{G}_h.
\end{align*}

In order to construct   Bayes and MAP tests, we assume that for given $h\in \mathcal{H}$ 
\begin{itemize}
\item  the set $\{\theta_{h,t},\ t\in \mathcal{G}_h\}$ contains the only one negative  
entry
$\theta_{h,\tau}$;
\item  $\tau$ is an unobservable random variable uniformly distributed on $\mathcal{G}_h$.
\end{itemize}

%Since  $\theta_{h,\tau}$  is unknown, we make use of an  improper prior  for this parameter.
\subsection{A Bayes test}
With the arguments used in deriving  \eqref{eq5}, we get  the  following  Bayes test:   
 $\mathbf{H}_0^h$ is rejected  if
\begin{equation*}
\begin{split}
 \frac{1}{n_h}\sum_{t\in \mathcal{G}_h}S\biggl(\frac{\hat{\theta}_{h,t}}{\sigma_h}\biggr)
\ge t_\alpha^B,
\end{split}
\end{equation*}
where  $S(\cdot)$  is defined by \eqref{eq6}.
The critical level $t_\alpha^B$ is defined by a conservative way, i.e., as a solution to 
\[
\max_{\Theta\ge 0}\mathbf{P}_\Theta\biggl\{ \frac{1}{ n_h} \sum_{t\in \mathcal{G}_h} S\biggl(\frac{\hat{\theta}_{h,t}}{\sigma_h}\biggr) > t_\alpha^B\biggr\}=\alpha,
\]
where here
$\mathbf{P}_\Theta $ stands for the measure generated by observations $\hat{\theta}_{h,t}$ defined by \eqref{eq2} for 
given $\Theta=\{\theta_{h,t}, h\in \mathcal{H}, t\in \mathcal{G}_h\}$. 

It follows  from Mudholkar's  theorem \cite{M-1966}, see also Theorem  6.2.1 in \cite{T}, that for any $\Theta$  with nonnegative entries $\theta_{h,t}\ge 0$
\begin{equation}\label{eq9}
\mathbf{P}_\Theta\biggl\{ \frac{1}{ n_h} \sum_{t\in \mathcal{G}_h} S\biggl(\frac{\hat{\theta}_{h,t}}{\sigma_h}\biggr) > x \biggr\} \le  \mathbf{P}\biggl\{\frac{1}{ n_h}\sum_{t\in \mathcal{G}_h} S(\xi_{h,t})\ge x\biggr\}
\end{equation}
and, thus, $t_\alpha^B$ may be computed as a solution to
\begin{equation}\label{eq10}
\mathbf{P}\biggl\{\frac{1}{ n_h}\sum_{t\in \mathcal{G}_h} S(\xi_{h,t})\ge t_\alpha^B\biggr\}
=\alpha.
\end{equation}
Therefore our next step is to study the following random variable:
\[
B_h(\xi) \stackrel{\rm def}{=} 
\frac{1}{ n_h} \sum_{t\in \mathcal{G}_h} S(\xi_{h,t}).
\] 
\subsubsection{A weak approximation of $B_h(\xi)$}
We begin with computing a weak limit of $B_h(\xi)$ as $h\rightarrow0$. 
Recall some standard definitions (see, e.g., \cite{No}).
\begin{definition}
Let $X_1$ and $X_2$ be independent copies of a random variable $X$. Then $X$ is said to be stable if for any constants $a > 0$ and $b > 0$ the random variable $aX_1 + bX_2$ has the same distribution as $cX + d$ 
for some constants $c > 0$ and $d$.
\end{definition}
 In the class of stable distributions there is an interesting sub-class of the so-called stable distributions with the index of stability $\alpha=1$. For brevity,
 we will call them \textit{1-stable distributions}. The formal definition of this class is as follows:  
\begin{definition}
A random variable $X$ is called 1-stable if its characteristic function can be written as
\begin{equation}\label{eq11}
\mathbf{E}\exp({\rm  i}tX)=
\exp\biggl(\mu {\rm i} t-|ct|-{\rm i}\frac{2\beta |c|}{\pi} t\log(|t|)
\biggr).
\end{equation}
\end{definition}

The next theorem shows that the weak limit of $B_h(\xi)-\log(n_h)$ is a 1-stable distribution.
\begin{theorem} \label{th1}
\begin{equation*}
\begin{split}
& \lim_{h\rightarrow0} \mathbf{E}\exp\bigl\{{\rm i}t
\bigl[B_h(\xi)-\log(n_h)+\gamma\bigr]\bigr\}
 =
\exp\biggl\{ {\rm i}t\log\frac{1}{|t|}
-\frac{\pi|t|}{2}\biggr\}, 
\end{split}
\end{equation*}
where $\gamma\approx 0.57721$ is Euler's constant.
\end{theorem}

In other words, this theorem states that 
\begin{equation*}
\lim_{h\rightarrow0}\bigl[B_h(\xi)-\log(n_h)+\gamma]\stackrel{\mathcal{D}}{=}\zeta,
\end{equation*}
where $\zeta$ is a  1-stable random variable (see \eqref{eq11}) with 
\begin{equation}\label{eq.12}
\mu=0, \  c=\frac{\pi}{2},\  \beta =1.
\end{equation}
Apparently, $\zeta$ appeared firstly in \cite{D}. Emphasize also that this random variable  originate usually in  Bayes hypothesis testing related to sparse vectors, see e.g. \cite{BB}, \cite{IS}.

The probability distribution of $\zeta$ has the following invariance property  
that plays  an important role  in Bayes tests aggregation. 
\begin{proposition} Let  $\zeta_k $ be i.i.d. copies of $\zeta$ and $\bar{\pi}$ be a probability distribution on 
$\mathbb{Z}^+$ with a bounded entropy. Then
\begin{equation}\label{eq13}
\begin{split}
%&\mathbf{E}\exp(it\zeta)=\exp\biggl\{ {\rm i} t +{\rm i}t\log\frac{1}{|t|}
%-\frac{\pi|t|}{2}\biggr\},\\
&\sum_{k=1}^\infty\bar{\pi}_k\biggl(\zeta_k-\log\frac{1}{\bar{\pi}_k}\biggr) \stackrel{\mathcal{D}}{=}  \zeta.
\end{split}
\end{equation}
\end{proposition}
The proof of \eqref{eq13} follows immediately from \eqref{eq11} and \eqref{eq.12}. 

\subsubsection{A strong approximation of $B_h(\xi)$}
Theorem \ref{th1} is not very informative about the tail behavior of the distribution of $B_h(\xi)$. However, for obtaining a good  approximation of $t_\alpha^B$ in \eqref{eq10} this behavior may play a crucial role
because in some applications $\alpha$ may be very small (of order $10^{-7}$) and so,  the Monte-Carlo method  and Theorem \ref{th1} may not be good in this case.

 Therefore our goal is to find
 an approximation of $
 B_h(\xi)
$ that
 controls well the tail of its distribution. 
Fortunately, this can be easily done. 
   It is clear that
   \[
   \Phi(\xi_k)\stackrel{\mathcal{D}}{=}U_k,
   \]
   where $U_k$ are i.i.d. random variables uniformly distributed on $[0,1]$. Hence
   \[
   B_h(\xi)\stackrel{\mathcal{D}}{=}\frac{1}{n_h}\sum_{k=1}^{n_h} \biggl[\frac{1}{U_k}-1\biggr]=
   \frac{1}{n_h}\sum_{i=1}^{n_h} \frac{1}{U_{(k)}}-1,
   \]
  where $U_{(k)}$ is a non-decreasing permutation of $U_k, \ k=1,\ldots,n_h$.
   The distribution of 
$U_{(1)}, U_{(2)},\ldots, U_{(n_h)}$ can be easily  obtained  with  the help of the Pyke theorem \cite{P}
\begin{equation}\label{eq14}
U_{(k)}\stackrel{\mathcal{D}}{=}\frac{\mathcal{E}_{k}}{\mathcal{E}_{n_h+1}},
\end{equation}
\text{where} 
\begin{equation*}
\mathcal{E}_{k}=\sum_{l=1}^k \varkappa_l
\end{equation*}
is the cumulative sum of i.i.d. standard exponentially distributed random variables $\varkappa_l$
\[
\mathbf{P}\bigl\{\varkappa_l\ge y\bigr\}=\exp(-y).
\]
In other words, $\mathcal{E}_k \sim {\rm Gamma}(k,1)$.
With this in mind, we obtain
    \begin{equation}\label{eq15}
    \begin{split}
   B_h(\xi)\stackrel{\mathcal{D}}{=}& \biggl[1+O\biggl(\frac{1}{\sqrt{n_h}}\biggr)\biggr]
   \sum_{k=1}^{n_h} \frac{1}{\mathcal{E}_k}-1  \\=&
    \biggl[1+O\biggl(\frac{1}{\sqrt{n_h}}\biggr)\biggr]\biggl[
   \sum_{k=1}^{n_h} \biggl(\frac{1}{\mathcal{E}_k}-\frac{1}{k}\biggr) + \sum_{k=1}^{n_h} \frac{1}{k}\biggr]-1.
  \end{split}
   \end{equation}
   
  Next, we make use of the following simple equations:
   \[
\biggl\{ \mathbf{E}  \biggl[  \sum_{k=n_h+1}^\infty \biggl(\frac{1}{\mathcal{E}_k}-\frac{1}{k}\biggr)\biggr]^{2m}\biggr\}^{1/(2m)}\le O\biggl(\frac{1}{\sqrt{n_h}}\biggr)
   \]
   and 
   \[
   \sum_{k=1}^{n_h} \frac{1}{k}=\log(n_h)+\gamma +O\biggl(\frac{1}{n_h}\biggr).
   \]
So, substituting them in \eqref{eq15}, we arrive at  the following theorem. 
\begin{theorem}\label{th3}  Let
\begin{equation}\label{eq16}
\zeta^\circ =\sum_{k=1}^\infty\biggl(\frac{1}{\mathcal{E}_k}
-\frac{1}{k}\biggr).
\end{equation}
Then
\begin{equation}\label{eq17}
B_h(\xi)-\log(n_h)+\gamma \stackrel{\mathcal{D}}{=}  \bigl[1+O(\varepsilon_h)\bigr]\zeta^\circ
+2\gamma-1 + O( \varepsilon_h)\log(n_h),
\end{equation}
where $\varepsilon_h$ is such that
\begin{equation}\label{rate}
\bigl[\mathbf{E}\bigl(\varepsilon_h^{2m}\bigr)\bigr]^{1/(2m)}\le \frac{C_m}{\sqrt{n_h}}.
\end{equation}
\end{theorem}

\begin{rem}
 The  random variable $\zeta$ in Theorem \ref{th1}
admits the following representation
\[
\zeta=\sum_{k=1}^\infty\biggl(\frac{1}{\mathcal{E}_k}
-\frac{1}{k}\biggr)+2\gamma-1.
\]
 Notice also that it follows immediately from  \eqref{eq17} that  convergence rate in  Theorem \ref{th1} is  $\log(n_h)/\sqrt{n_h}$, i.e., 
 as $h\rightarrow0$,
 \begin{equation*}
\begin{split}
&  \mathbf{E}\exp\bigl\{{\rm i}t
\bigl[B_h(\xi)-\log(n_h)+\gamma\bigr]\bigr\}
 =
\exp\biggl\{ {\rm i}t\log\frac{1}{|t|}
-\frac{\pi|t|}{2}+O\biggl(\frac{\log(n_h)}{\sqrt{n_h}}\biggr)\biggr\}.
\end{split}
\end{equation*}
\end{rem}

\begin{comment}
 \begin{definition}  Non-negative random sequences
 $X_n$ and  $Y_n$  are called asymptotically  $\log$-equivalent 
 $
 X_n\stackrel{\log}{\sim} Y_n
 $
 if 
 \[
\lim_{n\rightarrow\infty}\sup_{x\ge 0}\Bigl|\log \mathbf{P}\{X_n\ge x\}-\log \mathbf{P}\{Y_n\ge x\}\Bigr|=0.
 \]
 \end{definition}
\end{comment}

\begin{figure}
\includegraphics[angle=0,width=0.49\textwidth,height=0.25\textheight]{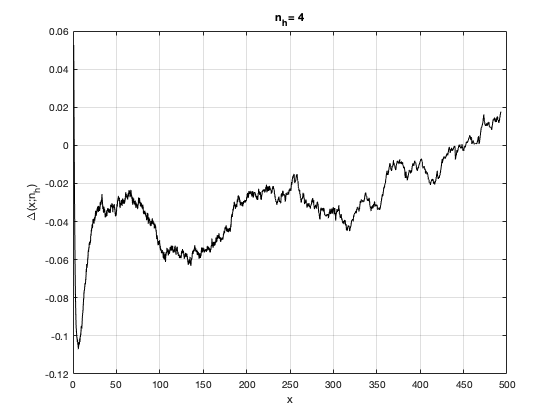}
\hfill
\includegraphics[angle=0,width=0.49\textwidth,height=0.25\textheight]{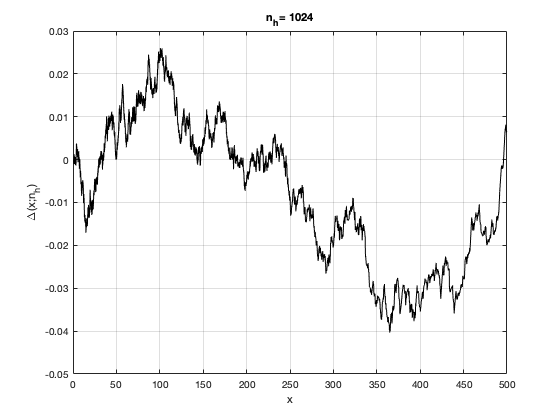}
 \caption{Log-tail approximation errors  $\Delta(x;n_h)$ for $n_h=4$ and $n_h=1024$.\label{F.1}}
\end{figure}

Figure \ref{F.1} illustrates   numerically Theorem \ref{th3} and the above remark showing log-tail approximation error $$
\Delta(x;n_h) =\log\Bigr[\mathbf{P}\bigl\{B_h\{\xi) -\log(n_h)+\gamma\ge x\bigr\}\Bigr]-\log\Bigr[\mathbf{P}\bigl\{\zeta^\circ+2\gamma-1\ge x\bigr\}\Bigr].
 $$ 
 computed with the help of the Monte-Carlo method with $0.5 \cdot 10^{6}$ replications.
This picture shows  that even for small $n_h=4$ the approximation \eqref{eq17} works very good.

\subsection{ A MAP test} 
Similarly to the  Bayes test, we can construct  the MAP test that   rejects $\mathbf{H}_0^h$
if
\begin{equation*}
\max_{t\in \mathcal{G}_h}\frac{1}{n_h}S\biggl(\frac{\hat{\theta}_{h,t}}{\sigma_h}\biggr)\ge t_\alpha^M, 
\end{equation*}
where $t_\alpha^M$ is defined as a solution to
\[
\max_{\Theta\ge 0}\mathbf{P}_\Theta\biggl\{  \max_{t\in \mathcal{G}_h}\frac{1}{n_h}S\biggl(\frac{\hat{\theta}_{h,t}}{\sigma_h}\biggr)> t_\alpha^M\biggr\}=\alpha.
\]
Similarly to \eqref{eq9},  $t_\alpha^M$ may be obtained from
\[
\mathbf{P}\biggl\{  \max_{t\in \mathcal{G}_h}\frac{S(\xi_{h,t})}{n_h}> t_\alpha^M\biggr\}=\alpha.
\] 
  
   As to  the limit distribution of $\max_{t\in \mathcal{G}_h}{S(\xi_{h,t})}/{n_h}$, as $h\rightarrow0$, it follows 
   immediately from \eqref{eq14} that
 \begin{equation}\label{eq18}
  \max_{t\in \mathcal{G}_h}\frac {S(\xi_{h,t})}{n_h}\stackrel{\mathcal{D}}{=} \frac{1+\varepsilon_h}{\varkappa},\ \text{as}\ h\rightarrow0,
  \end{equation}
  where  $\varkappa$ is a standard exponential random variable and  $\varepsilon_h$  satisfies
  \eqref{rate}.

\section{Multi-level testing}
\subsection{MAP multi-level tests}

A heuristic idea behind   our construction of multi-level MAP tests for \eqref{eq3}  is related to  \eqref{eq18} and consists in  
computing a positive deterministic function $U_h,\, h\in\mathcal{H},$ 
bounding from above the random  process $\log(1/\varkappa_h),\,  h\in \mathcal{H},$
where $\varkappa_h$ are independent standard exponential random variables. 
In other words, we are looking for $U_h$
such that 
\[
\zeta^U=\sup_{h\in\mathcal{H}}\biggl[\log\frac{1}{\varkappa_h}-U_h\biggr]
\]
would be a non-degenerate random variable.

Let $q_\alpha^U$ be $\alpha$-value of  $\zeta^U$,
i.e.,  solution to
\[
\mathbf{P}\bigl\{\zeta^U\ge q_\alpha^U\bigr\}=\alpha.
\]

Therefore with \eqref{eq18}, upper bounding  random process $\log(1/\varkappa_h)$  by $U_h$, we arrive at 
the   test that rejects
  $\mathbf{H}_0$ if
\begin{equation}\label{eq19}
\sup_{h\in\mathcal{H}}\biggl\{\max_{t\in \mathcal{G}_h}\log\biggl[\frac{1}{n_h}S\biggl(\frac{\hat{\theta}_{h,t}}{\sigma_h}\biggr)\biggr] -U_h\biggr\}\ge q_\alpha^U.
\end{equation}

Computing $q_\alpha^U$ is based  on the following simple fact. Assume that
\[
K^U=\log\biggl[\sum_{h\in\mathcal{H}}{\rm e}^{-U_h}\biggr]<\infty.
\]
Then
\begin{equation}\label{eq20}
\sup_{h\in\mathcal{H}}\biggl[\log\frac{1}{\varkappa_h}-U_h\biggr]-K^U\stackrel{\mathcal{D}}{=}
\log \frac{1}{\varkappa} .
\end{equation}
The proof of this identity  is very simple. Indeed,
\begin{equation*}\label{P2.p1}
\begin{split}
&\mathbf{P}\Bigl\{\sup_{h\in\mathcal{H}}\bigl[\log\frac{1}{\varkappa_h}-U_h\bigr]-K^U> x \Bigr\}
\\&\quad =
1-\prod_{h\in \mathcal{H}}\mathbf{P}\Bigl\{\log\frac{1}{\varkappa_h}\le U_h+x+K^U\Bigr\}
\\ &\quad =1-\exp\biggl\{-\sum_{h\in \mathcal{H}}\exp\biggl[-x -U_h-K^U\biggr]\biggr\}
 =1-\exp[-\exp(-x)].
\end{split}
\end{equation*}

Let us we denote
\begin{equation*}
\bar{\pi}_h={\rm e}^{-U_h} \biggl/ \sum_{h\in \mathcal{H}}{\rm e}^{-U_h} ,
\end{equation*}
then \eqref{eq20}  can be rewritten in the following form:
\begin{proposition}\label{P2}  Let $\bar{\pi}$ be a probability distribution on $\mathcal{H}$. Then
\begin{equation}\label{eq21}
\sup_{h\in\mathcal{H}}\biggl[\log\frac{1}{\varkappa_h}+\log(\bar{\pi}_h)\biggr]\stackrel{\mathcal{D}}{=}
\log \frac{1}{\varkappa}.
\end{equation}
\end{proposition}

Therefore with the help of \eqref{eq20} we can  compute  $\alpha$-critical level $q_\alpha^U$ in \eqref{eq19}
\[
q_\alpha^U=q_\alpha^\varkappa+K^U,
\]
where
\[
q_\alpha^\varkappa=-\log\biggl(\log\frac{1}{1-\alpha}\biggr)
\]
is $\alpha$-value of $\log(1/\varkappa)$. 

Summarizing (see \eqref{eq19}),   the  MAP multi-level test rejects $\mathbf{H}_0$ if
\begin{equation}\label{eq22}
\sup_{h\in\mathcal{H}}\bigl\{Z^M_h+\log(\bar{\pi}_h)\bigr\}\ge q_\alpha^\varkappa,
\end{equation} 
where
\[
Z_h^M=\max_{t\in \mathcal{G}_h}\log\biggl[\frac{1}{n_h}S\biggl(\frac{\hat{\theta}_{h,t}}{\sigma_h}\biggr)\biggr]
\]
and $\bar{\pi}$ is a probability distribution on $\mathcal{H}$.

In order to study the performance of this method, we analyze  the type II error probability. 
 For given $\{\rho,\tau: \rho \in\mathcal{H}, \ \tau\in \mathcal{G}_g\}$ 
  and $A \in\mathbb{R}^+$ define    
 \begin{equation}\label{Q}
 \Theta_{\rho,\tau}(A)=\Bigl\{\theta_{h,t}: \theta_{\rho,\tau}=-A; \ \theta_{h,t}\ge 0, (h,t)\neq (\rho,\tau)\Bigr\}.
 \end{equation}
 In other words, we consider 
 the situation, where all shifts $\theta_{h,t}$ in \eqref{eq2}  are positive except the only one.
 The position  of the negative entry  $\{\rho, \tau\}$ and its amplitude are unknown, but it is assumed that  $\{\rho, \tau\}$ are  random variables with the
distribution defined by
\begin{itemize}
\item $\mathbf{P}\{\rho =h\}=\bar{\pi}_h$,
\item $\mathbf{P}\{\tau=t|\rho=h\}=n_h^{-1}$,
\end{itemize} 
   where $\bar{\pi}$ is a probability distribution on $\mathcal{H}$ with a bounded entropy
   \[
 H_{\bar{\pi}}=\sum_{h\in\mathcal{H}}\bar{\pi}_h\log \frac{1}{\bar{\pi}_h}.
 \]

 In what follows, we will  deal with priors  $\bar{\pi}$  with large  uncertainties assuming 
 that $\bar{\pi} \rightarrow 0$, or more precisely,
 $
 \sup_{h\in\mathcal{H}}\bar{\pi}_h \rightarrow 0,
 $
but such that
 \begin{equation}\label{eq23}
 \lim_{\bar{\pi}\rightarrow0}\frac{1}{\log[H_{\bar{\pi}}]}\sum_{h\in\mathcal{H}}\bar{\pi}_h\biggr|H_{\bar{\pi}}-\log \frac{1}{\bar{\pi}_h}\biggr|=0.
  \end{equation}
  
 In particular, we  will consider  the following class of    prior distributions: 
 \begin{equation}\label{prior}
\bar{\pi}_h= \bar{\pi}_h^{\omega,\nu}= \nu \biggl[\frac{\log_2(1/h)}{\omega}\biggr]\bigg/
 \sum_{k=1}^\infty \nu\biggl(\frac{k}{\omega}\biggr)\approx \frac{1}{\omega}\nu\biggl[\frac{\log_2(1/h)}{\omega}\biggr].
 \end{equation}
 This class is characterized by the  bandwidth $\omega>1$ and the probability density
 $\nu(x),\ x\in \mathbb{R}^+$,  which is assumed to be  continuous, bounded, and with
 \begin{equation} \label{eq31}
H_{\nu}= \int_0^\infty \nu(x)\log \frac{1}{\nu(x)} dx<\infty,
 \quad
 \int_0^\infty \nu(x)\log(x+1)\,dx <\infty.
 \end{equation}
 A typical example of a such distribution is the uniform one that corresponds to $\nu(x)=1,\, x\in [0,1]$.
 
   It is clear that $\bar{\pi}_h^{\omega,\nu}\rightarrow0$ as   $\omega \rightarrow\infty$ and that Condition \eqref{eq23} holds. 
 
 Let us  begin with the case, where the prior distribution   is known, the case of unknown 
   $\bar{\pi}$ will be considered later in Section \ref{robust}.

    The  type II error probability over $ \Theta_{\rho,\tau}(A)$ of the MAP test  \eqref{eq22} is defined as follows:
 \[
 \beta^M_{\rho,\tau}(A)=\sup_{\Theta\in \Theta{\rho,\tau}(A)}\mathbf{P}_\Theta\Bigl\{ \max_{h\in \mathcal{H}}\bigl[Z_h^M +\log(\bar{\pi}_h)\bigr]
 \le q_\alpha^\varkappa\Bigr\}.
 \]
 Our goal is to study the average type II error probability
 \[
\bar{ \beta}^M_{\bar{\pi}}(\mathbf{A})=\sum_{h\in \mathcal{H}} \frac{\bar{\pi}_h}{n_h}\sum_{t\in \mathcal{G}_h}\beta^M_{h,t}(A_h),
 \]
  where here and below $\mathbf{A}=\{A_h, \, h\in \mathcal{H}\}$.
  
Denote for brevity
\[
R_h(q,H)=2[q+\log(n_h)+H]-\log[4\pi(q+\log(n_h)+H)]
\]
and  
\begin{equation*}
\begin{split}
\log^*(x)=&\log[\log(x)], \quad
H^*_{\bar{\pi}}=\log(H_{\bar{\pi}}).
\end{split}
\end{equation*}

The next theorem shows that  $R_h(q_\alpha^\varkappa,H_{\bar{\pi}})$ 
 is a critical signal/noise ratio.  
 Roughly speaking, this means that 
 if
 \[
 \frac{A_h}{\sigma_h} \stackrel{\bar{\pi}}{\le}\sqrt{ R_h(q_\alpha^\varkappa,H_{\bar{\pi}})}+x
 \]
 for any given $x>0$,
 then the MAP multi-level test cannot discriminate between $\mathbf{H}_0$ and    $\mathbf{H}_1$.
 Otherwise, if 
 \[
 \frac{A_h}{\sigma_h}  \stackrel{\bar{\pi}}{\ge }\sqrt{R_h(q_\alpha^\varkappa,H_{\bar{\pi}})}+\sqrt{\epsilon H^*_{\bar{\pi}}},
 \]
 for some $\epsilon>0$,
 then reliable testing is possible.  
 
  In the next theorem,  $ \mathbf{E}_{\bar{\pi}}$ stands for the expectation w.r.t. $\bar{\pi}$.
  \begin{theorem} \label{th6}
 Suppose \eqref{eq23} holds.
 If  for some $x\in \mathbb{R}$ and $\epsilon>0$
 \begin{equation}\label{eq24} 
\lim_{\bar{\pi}\rightarrow 0}\frac{1}{H^*_{\bar{\pi}}}  \mathbf{E}_{\bar{\pi}}\biggl[\biggl(\frac{A_{h}}{\sigma_h}-x\biggr)^2
+\epsilon H^*_{\bar{\pi}} -R_h(q_\alpha^\varkappa,H_{\bar{\pi}}) \biggr]_+ =0,
 \end{equation}
 then
 \begin{equation}\label{eq25}
 \lim_{\bar{\pi}\rightarrow 0}  \bar{\beta}^M_{\bar{\pi}}(\mathbf{A}) 
\ge  (1-\alpha)[1-\Phi(x)].
 \end{equation}
 If for some $\epsilon>0$
 \begin{equation}\label{eq26} 
\lim_{\bar{\pi}\rightarrow 0}\frac{1}{H^*_{\bar{\pi}}} \mathbf{E}_{\bar{\pi}}\biggl[
R_h(q_\alpha^\varkappa,H_{\bar{\pi}})+2\sqrt{\epsilon H^*_{\bar{\pi}}}\frac{A_h}{\sigma_h}-\frac{A_{h}^2}{\sigma_h^2}\biggr]_+ =0,
\end{equation}
 then
 \begin{equation}\label{eq27}
 \lim_{\bar{\pi}\rightarrow 0} \bar{\beta}^M_{\bar{\pi}}(\mathbf{A}) =0.
 \end{equation}
 \end{theorem}
 
\subsection{Multi-level Bayes  tests}
To construct these tests, let us
 consider the following statistics:
 \[
 Z_h^B= \frac{1}{n_h}\sum_{t\in \mathcal{G}_h}S\biggl(\frac{\hat{\theta}_{h,t}}{\sigma_h}\biggr)
-\log(n_h)-\gamma+1,\ h\in \mathcal{H}.
 \] 
 When all $\theta_{h,t}=0$,
in view of Theorem \ref{th3}, these random variables are  approximated by the family of  independent and identically distributed  random variables $\zeta^\circ_h, h\in \mathcal{H}$, defined by \eqref{eq16}. An important
property of this family  is provided by 
 \eqref{eq13}, which is  used  in our construction multi-level Bayes tests. 
More precisely, the multi-level Bayes test rejects $\mathbf{H}_0$ if
\begin{equation*}
\sum_{h\in \mathcal{H}}\bar{\pi}_h  \biggl(Z_h^B-\log\frac{1}{\bar{\pi}_h}\biggr)\ge q_\alpha^\circ,
\end{equation*} 
where $q_\alpha^\circ$ is $\alpha$-value of $\zeta^\circ$.

  The  type II error probability   over  $\Theta_{\rho,\tau}(A)$ (see \eqref{Q}) is defined by 
  \[
  \beta^B_{\rho,\tau}(A)=\sup_{\Theta\in \Theta_{\rho,\tau}(A_\rho)}\mathbf{P}_\Theta\biggl\{ 
 \sum_{h\in \mathcal{H}}\bar{\pi}_h  \biggl(Z_h^B-\log\frac{1}{\bar{\pi}_h}\biggr)\le q_\alpha^\circ \biggr\}
 \]
 and our goal is to analyze  the average  type II error probability
  \[
 \bar{\beta}^B_{\bar{\pi}}(\mathbf{A})=\sum_{h\in \mathcal{H}} \frac{\bar{\pi}_h}{n_h}\sum_{\tau\in \mathcal{G}_h} \beta^B_{h,\tau}(A_h).
 \]

 \begin{theorem} \label{th9}
 Suppose \eqref{eq23} holds and  for some $x\in\mathbb{R}$ and $\epsilon>0$
 \begin{equation}\label{eq28} 
 \lim_{\bar{\pi}\rightarrow0} \frac{1}{ H^*_{\bar{\pi}} } \mathbf{E}_{\bar{\pi}} \biggl[\biggl(
\frac{A_h}{\sigma_h}-x\biggr)^2+\epsilon H^*_{\bar{\pi}} -R_h[\log(q_\alpha^\circ),H_{\bar{\pi}}]
 \biggr]_+=0,
 \end{equation}
 then
 \begin{equation}\label{eq29}
 \lim_{\bar{\pi}\rightarrow0} \bar{\beta}^B_{\bar{\pi}}(\mathbf{A}) \ge 
 (1-\alpha)[1-\Phi(x)].
 \end{equation}
 If  for some $\epsilon>0$ 
 \begin{equation}\label{eq30} 
\lim_{\bar{\pi}\rightarrow 0}\frac{1}{H^*_{\bar{\pi}}} \mathbf{E}_{\bar{\pi}}\biggl[
R_h\bigl[\log(q_\alpha^\circ),H_{\bar{\pi}}\bigr]+2\sqrt{\epsilon H^*_{\bar{\pi}}}\frac{A_h}{\sigma_h}-\frac{A_{h}^2}{\sigma_h^2} \biggr]_+ =0,
 \end{equation}
 then
 \begin{equation}\label{eq35.x}
 \lim_{\bar{\pi}\rightarrow0}\bar{ \beta}^B_{\bar{\pi}}(\mathbf{A})
  =0.
 \end{equation}
 \end{theorem}

\medskip
\noindent
\begin{rem}
Notice that as $\alpha\rightarrow0$
\[
\log(q_\alpha^\circ)=(1+o(1)) q_\alpha^\varkappa =(1+o(1)) \log\frac{1}{\alpha}.
\]
Therefore,  since $H_{\bar{\pi}}\rightarrow \infty$ as $\bar{\pi}\rightarrow0$, conditions \eqref{eq24}  and \eqref{eq28}  
along with \eqref{eq26} and  \eqref{eq30}   are almost equivalent. This means that in the considered statistical problem there is no substantial difference between 
MAP and Bayes tests. 
  \end{rem}

 \section{Adaptive  multi-level tests \label{robust}}
 The main drawback of the MAP and Bayes tests is related to their dependence on the  prior distribution $\bar{\pi}$ that is hardly known in practice.
 Therefore our next goal is to construct a test that, on the one hand,  does not depend on $\bar{\pi}$, but on the other hand, has a nearly 
 optimal critical signal-noise ratio.
 
 In order to simplify our presentation, we will deal with the class  of    prior distributions 
 $\bar{\pi}^{\omega,\nu}$ defined by \eqref{prior}.
The  entropy of $\bar{\pi}^{\omega,\nu}$   obviously satisfies
 \begin{equation}\label{eq32}
 H_{\bar{\pi}^{\omega,\nu}}=\log(\omega)+H_{\nu}+o(1), \ \omega\rightarrow\infty,
 \end{equation}
and therefore denote for brevity
\begin{equation}\label{eq33}
\widetilde{R}_h(q,\omega)=2[q+\log(n_h)+\log(\omega)]-\log[4\pi(q+\log(n_h)+\log(\omega)].
\end{equation}
\begin{comment}
Let 
\[
\beta^M(A;\bar{\pi}^{W,\nu})=\sum_{h\in \mathcal{H}}\frac{\bar{\pi}_h^{\omega,\nu}}{n_h}\sum_{t\in \mathcal{G}_h} \beta^M_{h,t}(A_h)
\]
be the average type II error probability.
\end{comment}

 With \eqref{eq32}, Condition \eqref{eq23} is checked  easily and 
the next result  follows immediately from Theorem \ref{th6}.
 \begin{cor}  \label{cor1} If for some $x\in\mathbb{R}$ and $\epsilon>0$
 \begin{equation*}\label{u.2}
 \begin{split}
\lim_{\omega\rightarrow\infty} \frac{1}{\log^*(\omega)} \mathbf{E}_{\bar{\pi}^{\omega,\nu}}
\biggl[\biggl(\frac{A_{h}}{\sigma_h}-x\biggr)^2+\epsilon \log^*(\omega)-\widetilde{R}_h(q_\alpha^\varkappa,\omega)\biggr]_+&=0,
\end{split}
 \end{equation*}
 then
 \[
\lim_{\omega\rightarrow\infty} \bar{ \beta}^M_{\bar{\pi}^{\omega,\nu}}(\mathbf{A}) \ge (1-\alpha)[1-\Phi(x)].
 \]
 If for some $\epsilon>0$
 \begin{equation*}
 \begin{split}
\lim_{\omega\rightarrow \infty}\frac{1}{\log^*(\omega)]} \mathbf{E}_{\bar{\pi}^{\omega,\nu}}\biggl[
\widetilde{R}_h(q_\alpha^\varkappa,\omega)+2\sqrt{\epsilon \log^*(\omega)}\frac{A_{h}}{\sigma_h}-\frac{A_{h}^2}{\sigma_h^2}\biggr]_+ =0,
\end{split}
 \end{equation*}
 then
 \begin{equation*}
 \lim_{\omega\rightarrow \infty} \bar{\beta}^M_{\bar{\pi}^{\omega,\nu}}(\mathbf{A})=0.
 \end{equation*}
 \end{cor}

In order to construct an adaptive  test, let us compute a nearly minimal function $U_h$ in \eqref{eq20}.  We begin with
\[
\psi_0(x)=1+\log(x), \quad x\in \mathbb{R}^+,
\]
and then  iterate this function $m$ times
\[
\psi_{l}(x)=\psi_0\bigl[\psi_{l-1}(x)\bigr],\  l=1,\ldots, m.
\]
Finally, for given  $\varepsilon\in(0,1)$,
define
\begin{equation}\label{eq34}
L^{m,\varepsilon}(k)=-\log\bigg\{
\frac{1}{\varepsilon[\psi_m(k)]^\varepsilon}-\frac{1}{\varepsilon[\psi_m(k+1)]^\varepsilon}\biggr\},\ k\in \mathbb{Z}^+.
\end{equation}
Since $\psi_{m}(1)=1$, it is clear that 
\[
\sum_{k=1}^\infty\exp[-L^{m,\varepsilon}(k)]= \frac{1}{\varepsilon}.
\]

\begin{figure}
\includegraphics[angle=0,width=0.49\textwidth,height=0.25\textheight]{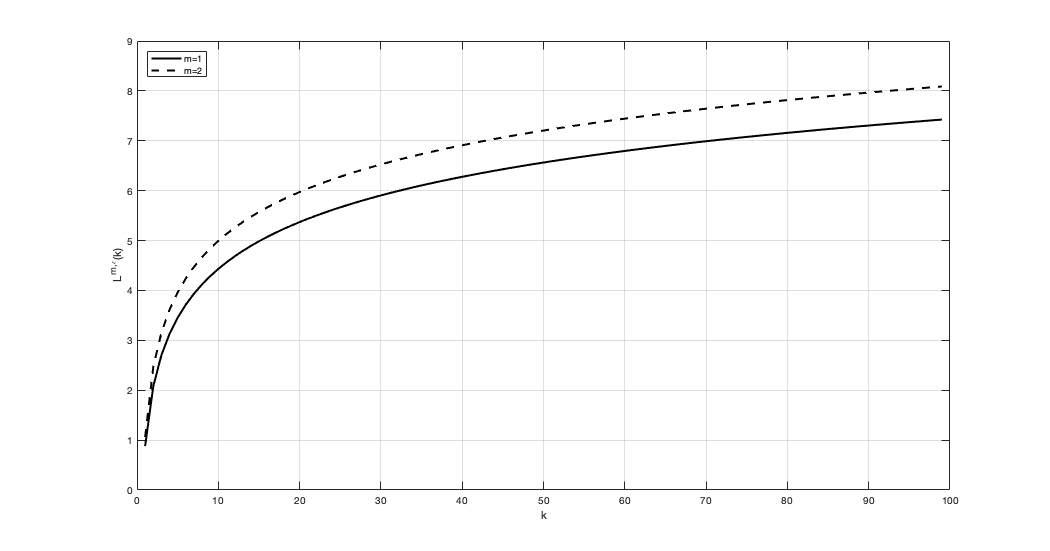}
\hfill
\includegraphics[angle=0,width=0.49\textwidth,height=0.25\textheight]{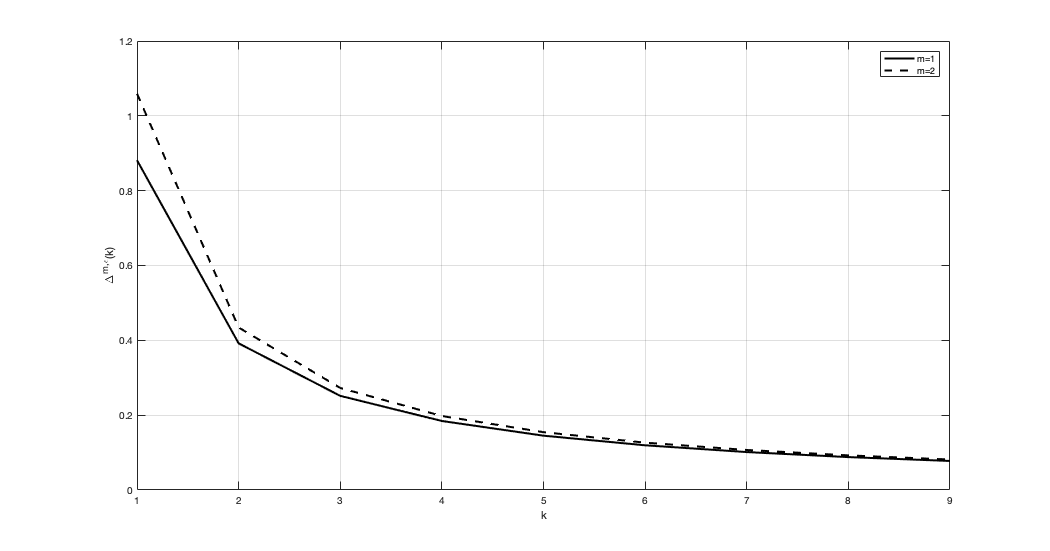}
 \caption{The functions $L^{m,\varepsilon}(\cdot)$   and  the approximation errors $\Delta^{m,\varepsilon}(\cdot)$ for  $m=1,2$ and $\varepsilon=0.1$.  \label{F.2}}
\end{figure}

In what follows, we will make use of  the following approximation of 
  $L^{m,\varepsilon}(k)$  for large $k$.
Denote (see \eqref{eq34})
 \begin{equation}\label{eq34.1}
\begin{split}
&\widetilde{L}^{m,\varepsilon}(k)=-\log\bigg[-\frac{1}{\varepsilon} \frac{[d\psi_m(k)]^{-\varepsilon}}{dk}
\biggr]=-\log\bigg[\frac{1}{[\psi_m(k)]^{1+\varepsilon}}\frac{d\psi_m(k)}{dk}
\biggr]\\&\qquad =\log(k)+\log[\psi_0(k)]+\cdots +\log[\psi_{m-1}(k)]+(1+\varepsilon)\log[\psi_m(k)].
\end{split}
\end{equation}
Since 
\[
\frac{d^2\psi_m(k)}{dk^2}\bigg/ \frac{d\psi_m(k)}{dk}=O\biggl(\frac{1}{k}\biggr),
\]
by the Taylor formula we obtain from \eqref{eq34}
\begin{equation}\label{eq34.2}
\Delta^{m,\varepsilon}(k)=\widetilde{L}^{m,\varepsilon}(k)-{L}^{m,\varepsilon}(k)=O\biggl(\frac{1}{\psi_m(k)}\biggr).
\end{equation}
Figure \ref{F.2} shows that  ${L}^{m,\varepsilon}(k)$ and  approximation errors $\Delta^{m,\varepsilon}(k)$.

Since $h\in \mathcal{H}=\{2^{-1},\ldots, 2^{-k},\ldots\}$, we  choose 
\begin{equation*}
U_h=L^{m,\varepsilon}[\log_2(1/h)]
\end{equation*}
and in view of \eqref{eq34} we arrive at the following prior distribution:
\begin{equation*}
{\Pi}_h=
\frac{1}{\{\psi_m[\log_2(1/h)]\}^\varepsilon}-\frac{1}{\{\psi_m[\log_2(1/h)+1]\}^\varepsilon}, \ h\in \mathcal{H}.
\end{equation*}

It is easy to check that the entropy of ${\Pi}$ is unbounded and this is why this distribution might be viewed as an improper  prior.

 The MAP  test  associated with
 $ {\Pi}$  rejects $\mathbf{H}_0$  when
\begin{equation*}
\max_{h\in \mathcal{H}}\bigl\{Z_h^M+\log\bigl({\Pi}_h\bigr)\bigr\}
 \ge q_\alpha^\varkappa
\end{equation*}
and its type II error probability over $ \Theta_{\rho,\tau}(A)$  (see \eqref{Q}) is defined by 
\[
\beta_{\rho,\tau}^{\Pi}(A)=\sup_{\Theta\in \Theta_{\rho,\tau}(A)}\mathbf{P}_\Theta\Bigl\{ \max_{h\in \mathcal{H}}\bigl[Z_h^M +\log({\Pi}_h)\bigr]
 \le q_\alpha^\varkappa\Bigr\}.
\]
Denote for brevity
\begin{equation}\label{R+}
R^+(q_\alpha^\varkappa,\omega)=\widetilde{R}(q_\alpha^\varkappa,\omega)+\log^*(\omega),
\end{equation}
where $\widetilde{R}(q_\alpha,\omega)$ is defined by \eqref{eq33}, and let 
\[
\bar{\beta}^{\Pi}_{\bar{\pi}^{\omega,\nu}}\bigl(\mathbf{A})=\sum_{h\in \mathcal{H}}\frac{\bar{\pi}_h^{\omega,\nu}}{n_h}\sum_{t\in \mathcal{G}_h} \beta^{\Pi}_{h,t}(A_h)
\]
be the average type II error probability.
 \begin{theorem}\label{th7}
 If for some $x\in\mathbb{R}$ and $\epsilon>0$
 \begin{equation*}
 \begin{split}
\lim_{\omega\rightarrow\infty} \frac{1}{\log^*(\omega)} \mathbf{E}_{\bar{\pi}^{\omega,\nu}}
\biggl[\biggl(\frac{A_{h}}{\sigma_h}-x\biggr)^2+\epsilon \log^*(\omega) -R_h^+(q_\alpha^\varkappa,\omega)\biggr]_+&=0,
\end{split}
 \end{equation*}
 then
 \[
\lim_{\omega\rightarrow\infty} \bar{\beta}^{\Pi}_{\bar{\pi}^{\omega,\nu}}(\mathbf{A})\ge (1-\alpha)[1-\Phi(x)].
 \]
 If for some $\epsilon>0$
 \begin{equation*} 
 \begin{split}
\lim_{\omega\rightarrow \infty}\frac{1}{\log^*(\omega)} \mathbf{E}_{\bar{\pi}^{\omega,\nu}}\biggl[
R_h^+(q_\alpha^\varkappa,\omega)+2\sqrt{\epsilon \log^*(\omega)}\frac{A_{h}}{\sigma_h} -\frac{A_{h}^2}{\sigma_h^2}\biggr]_+ =0,
 \end{split}
 \end{equation*}
 then
 \begin{equation*}
\lim_{\omega\rightarrow\infty}  \bar{\beta}^{\Pi}_{\bar{\pi}^{\omega,\nu}}(\mathbf{A})=0.
 \end{equation*}
 \end{theorem}

This theorem  and Corollary \ref{cor1} demonstrate  that the critical signal-noise ratio of the adaptive test 
 is  only slightly greater, see \eqref{R+},  (by the additive term $\log^*(\omega)$)  than the one of the MAP test that knows  
 the prior distribution $\bar{\pi}^{\omega,\nu}$.

\section{Appendix}

\subsection{Proof of Theorem \ref{th1}}

With a simple algebra we obtain
\begin{equation}\label{eq35}
\begin{split}
&\log\bigl\{ \mathbf{E}\exp\bigl[{\rm i}tB_h(\xi)\bigr]\bigr\} = n_h\log\biggl\{\frac{1}{\sqrt{2\pi}}\int_{-\infty}^\infty \cos
\biggl[\frac{tS(x)}{n_h }\biggr]
{\rm e}^{-x^2/2}\, dx\\
 &\qquad  +\frac{\rm i}{\sqrt{2\pi}}\int_{-\infty}^\infty\sin\biggl[\frac{tS(x) }{n_h}\biggr]{\rm e}^{-{x^2}/{2}}\, dx \biggr\}\\
&\ = n_h\log\biggl\{1+\int_{-\infty}^\infty \biggr[\cos
\biggl[\frac{t}{n_h }\biggl(\frac{1}{\Phi(x)}-1\biggr)\biggr]-1\biggl]
\, d\Phi(x)\\ &\qquad  +{\rm i}\int_{-\infty}^\infty\sin\biggl[\frac{t }{n_h}\biggl(\frac{1}{\Phi(x)}-1\biggr)\biggr]\, d\Phi(x) \biggr\}\\
&\ = n_h\log\biggl\{1+\int_{0}^1 \biggl\{\cos
\biggl[\frac{t}{n_h }\biggl(\frac{1}{u}-1\biggr)\biggr]-1\biggr\}
\, du \\ &\qquad  +{\rm i}\int_{0}^1\sin\biggl[\frac{t }{n_h}\biggl(\frac{1}{u}-1\biggr)\biggr]\, du \biggr\}\\ &\
= n_h\log\biggl\{1+\int_{0}^\infty\frac{1}{(1+u)^2}\biggl[ \cos
\biggl(\frac{tu}{n_h }\biggr)-1\biggr]
\, du  \\ & \qquad \qquad+{\rm i}\int_{0}^\infty
\frac{1}{(1+u)^2}\sin\biggl(\frac{tu }{n_h}\biggr)\, du \biggr\}\\
&\
= n_h\log\biggl\{1+\frac{|t|}{n_h}\int_{0}^\infty\frac{\cos(u)-1}{(|t|/n_h+u)^2}
\, du  +\frac{{\rm i}t}{n_h}\int_{0}^\infty
\frac{\sin(u)}{(|t|/n_h+u)^2}\, du \biggr\}.
\end{split}
\end{equation}
It is clear that  as $n_h\rightarrow \infty$
\begin{equation}\label{eq36}
\int_{0}^\infty\frac{\cos(u)-1}{(|t|/n_h+u)^2}
\, du \rightarrow 
\int_{0}^\infty\frac{\cos(u)-1}{u^2}
\, du = -\frac{\pi}{2}.
\end{equation}
Let   us choose $\epsilon_h<1$ such that 
$$
\lim_{h\rightarrow 0}\epsilon_h=0, \quad \lim_{h\rightarrow0}\epsilon_h n_h=\infty.
$$
Then we get by the  Taylor formula 
\begin{equation}\label{eq37}
\begin{split}
&\int_{0}^\infty
\frac{\sin(u)}{(|t|/n_h+u)^2}\, du=\int_{0}^{\epsilon_h}
\frac{\sin(u)}{(|t|/n_h+u)^2}\, du+
\int_{\epsilon_h}^\infty
\frac{\sin(u)}{(|t|/n_h+u)^2}\, du\\
&\qquad =  \int_{0}^{\epsilon_h}
\frac{u}{(|t|/n_h+u)^2}\, du+  O\biggl[\int_{0}^{\epsilon_h}
\frac{u^3}{(|t|/n_h+u)^2}\, du\biggr]\\ &\qquad\quad +(1+o(1))
\int_{\epsilon_h}^\infty
\frac{\sin(u)}{u^2}\, du\\
&\qquad =\log\bigg(1+\frac{\epsilon_h n_h}{|t|}\biggr)-\frac{\epsilon_h n_h}{(|t|+\epsilon_h n_h)}\\ & \qquad \quad+(1+o(1))
\int_{\epsilon_h}^\infty
\frac{\sin(u)}{u^2}\, du +O(\epsilon_h^2).
\end{split}
\end{equation}

 Next, integrating by parts, we obtain
 \begin{equation*}
 \begin{split}
  \int_{x}^\infty\frac{\sin(z)}{z^2}\, dz=&-\int_{x}^\infty \frac{\sin(z)}{z}\, d\biggl( \log\frac{1}{z}\biggr)
  \\=&\frac{\sin(x)}{x}\log \frac{1}{x}+\int_{x}^\infty \frac{z\cos(z)-\sin(z)}{z^2} \log \frac{1}{z}\, dz.
  \end{split}
 \end{equation*}
 Hence, as $x\rightarrow 0$
 \begin{equation*}
 \begin{split}
 \int_{x}^\infty\frac{\sin(z)}{z^2}\, dz=&\log\frac{1}{x}+\int_{0}^\infty \frac{z\cos(z)-\sin(z)}{z^2} \log \frac{1}{z}\, dz 
 +O\bigg(x^2\log \frac{1}{x}\biggr)\\
=&\log\frac{1}{x}+(1-\gamma)
 +O\bigg(x^2\log \frac{1}{x}\biggr),
 \end{split}
 \end{equation*}
 where $\gamma$ is Euler's constant.

With this equation we  continue \eqref{eq37} as follows:
\begin{equation*}
\begin{split}
\int_{0}^\infty
\frac{\sin(u)}{(|t|/n_h+u)^2}\, du=\log\frac{1}{|t|}+\log(n_h)-\gamma+O\biggl(\frac{|t|}{\epsilon_hn_h}\biggr)
+O\biggl(\epsilon_h^2\log\frac{1}{\epsilon_h}\biggr).
\end{split}
\end{equation*}

Substituting this equation  and \eqref{eq36} in \eqref{eq35}, we get
\begin{equation*}\label{c.3.a}
\begin{split}
\log\bigl\{ \mathbf{E}\exp\bigl[{\rm i}tB_h(\xi)\bigr]\bigr\} = -\frac{\pi |t|}{2}+
{\rm i}t\biggl(\log\frac{1}{|t|}+\log(n_h)-\gamma \biggr)+o(1),
\end{split}
\end{equation*}
thus, proving the theorem.

\subsection{Proof of Theorem \ref{th6}}
\textbf{I. A lower bound.} By \eqref{eq21} we  have for any given $x$ 
\begin{equation}\label{eq38}
\begin{split}
&\beta^M_{\rho,\tau}(A)  
\ge \mathbf{P}\biggl\{ \biggl[\log\biggl[  S\biggl(-\frac{A}{\sigma_\rho}+\xi_{\rho,\tau}\biggr)\biggr]-\log\frac{n_\rho}{\bar{\pi}_\rho}\biggr]\\&\qquad\bigvee
\biggl[\zeta_\rho-\log\frac{1}{\bar{\pi}_\rho}\biggr]\le q_\alpha^\varkappa\biggr\}
%\\
%&\quad \times
  \mathbf{P}\biggl\{\max_{h\in\mathcal{H}} \biggl[\zeta_h-\log\frac{1}{\bar{\pi}_h}\biggr] \le q_\alpha^\varkappa\biggr\}\\
&= \mathbf{P}\biggl\{   S\biggl(-\frac{A}{\sigma_\rho}+\xi_{\rho,\tau}\biggr)\le \exp\biggl[ q_\alpha^\varkappa
+\log\frac{n_\rho}{\bar{\pi}_\rho}\biggr]\biggr\}\\ 
&\quad \times 
 \mathbf{P}\biggl\{ \zeta_\rho \le q_\alpha^\varkappa+\log\frac{1}{\bar{\pi}_\rho}\biggr\}
 \mathbf{P}\bigl\{ \zeta \le q_\alpha^\varkappa \bigr\}
 \\
 &\ge (1-\alpha)^{1+\bar{\pi}_\rho} \mathbf{P}\biggl\{    S\biggl(-\frac{A}{\sigma_\rho}+\xi_{\rho,\tau}\biggr)\le 
 \exp\biggl[ q_\alpha^\varkappa+\log\frac{n_\rho}{\bar{\pi}_\rho}\biggr]; \xi_{\rho,\tau}\ge x \biggr\}\\
   &\ge (1-\alpha)^{1+\bar{\pi}_\rho} \mathbf{P}\bigl\{\xi_{\rho,\tau}\ge x \bigr\}
%  \\ &\qquad \times 
   \mathbf{1}\biggl\{ 
      S\biggl(-\frac{A}{\sigma_\rho}+x\biggr)\le\exp\biggl[ q_\alpha^\varkappa+\log\frac{n_\rho}{\bar{\pi}_\rho}\biggr]\biggr\}.
\end{split}
\end{equation} 

Let $R(z)\ge 0$ be a solution to
\[
S\Bigl[-\sqrt{R(z)}\Bigr]=z.
\]
It is easy to check with the help of \eqref{eq8} that as $z\rightarrow\infty$
\begin{equation}\label{eq39}
R(z)  =2\log(z)-\log[4\pi\log(z)]+o(1).
\end{equation}

Denote for brevity
\[
r_h(q,u)=  2\biggl(q+\log\frac{n_h}{u} \biggr) -
  \log\biggl[4\pi\biggl(q+\log\frac{n_h}{u}\biggr)\biggr].
\]

With \eqref{eq39} and the Markov inequality we obtain  for any $\epsilon>0$ 
\begin{equation} \label{eq40}
\begin{split}
&\mathbf{E}_{\bar{\pi}}\biggl[
\mathbf{1}\biggl\{ 
      S\biggl(-\frac{A_h}{\sigma_h}+x\biggr)\le\exp\biggl[ q_\alpha^\varkappa+\log\frac{n_h}{\bar{\pi}_h}\biggr]\biggr\}\biggr]
 \\ 
    &\quad= 1- \mathbf{E}_{\bar{\pi}} \biggl[\mathbf{1}\biggl\{ S\biggl(-\frac{A_h}{\sigma_h}+x\biggr) 
    > \exp\biggl[q_\alpha^\varkappa+\log\frac{n_h}{\bar{\pi}_h} \biggr]  \biggr\} \biggr]
     \\ 
    &\quad = 1- \mathbf{E}_{\bar{\pi}}\biggl[ \mathbf{1}\biggl\{ \biggl(\frac{A_{h}}{\sigma_h}-x\biggr) ^2
    -r_h(q_\alpha^\varkappa,\bar{\pi}_h)  >o(1) \biggr\}\biggr] \\    
  &\quad\ge 1-\frac{1}{\epsilon H^*_{\bar{\pi}}}\mathbf{E}_{\bar{\pi}} \biggl[ \biggl(
    \frac{A_h}{\sigma_h}-x\biggr)^2-r_h(q_\alpha^\varkappa,\bar{\pi}_h) 
    +\epsilon H^*_{\bar{\pi}}\biggr]_+ 
     .
\end{split}
\end{equation}

With  Condition \eqref{eq23}  and simple algebra it is easy to check that
\begin{equation*}
\lim_{\bar{\pi}\rightarrow0} \frac{1}{H^*_{\bar{\pi}}}\mathbf{E}_{\bar{\pi}} \biggl|r_h(q_\alpha^\varkappa,\bar{\pi}) 
    -R_h(q_\alpha^\varkappa,H_{\bar{\pi}})\biggr| =0.
\end{equation*}

Therefore \eqref{eq25} follows from \eqref{eq38}, \eqref{eq40} and the above equation. 

\bigskip
\noindent
\textbf{II. An upper bound.} Since $S(x)$ is a decreasing  function, we have obviously for any $x\ge 0$
\begin{equation*}
\begin{split}
 \bar{\beta}^M_{\bar{\pi}}(\mathbf{A})  
\le& \mathbf{E}_{\bar{\pi}}\biggl[ \mathbf{P}\biggl\{  S\biggl(-\frac{A_h}{\sigma_h}+\xi \biggr)
\le  \exp\biggl[ q_\alpha^\varkappa+\log\frac{n_h}{\bar{\pi}_h}\biggr] \biggr\}\biggr] \\
 \le &\mathbf{E}_{\bar{\pi}} \biggl[\mathbf{P}\biggl\{  S\biggl(-\frac{A_h}{\sigma_h}+\xi\biggr)
\le  \exp\biggl[ q_\alpha^\varkappa+\log\frac{n_h}{\bar{\pi}_h}\biggr]; \xi\le x \biggr\}\biggr]\\
 +&\mathbf{E}_{\bar{\pi}} \biggl[\mathbf{P}\biggl\{  S\biggl(-\frac{A_h}{\sigma_h}+\xi\biggr)
\le  \exp\biggl[ q_\alpha^\varkappa+\log\frac{n_h}{\bar{\pi}_h}\biggr]; \xi> x \biggr\}\biggr]\\
     \le& \mathbf{E}_{\bar{\pi}} \biggl[\mathbf{1}\biggl\{  S\biggl(-\frac{A_h}{\sigma_h}+x\biggr)
\le  \exp\biggl[ q_\alpha^\varkappa+\log\frac{n_h}{\bar{\pi}_h}\biggr] \biggr\} \biggr]
+ \mathbf{P}\bigl\{\xi>x\bigr\}.
%\\
%&\qquad +\sum_{\rho\in\mathcal{H}}\bar{\pi}_\rho \mathbf{P}\biggl\{  S\biggl(-\frac{A_\rho}{\sigma_\rho}+\xi_{\rho,\tau}\mathbf{1}\{\xi_{\rho,\tau}>x\}\biggr)
%\le  \exp\biggl[ q_\alpha+\log\frac{n_\rho}{\bar{\pi}_\rho}\biggr] \biggr\}.
 \end{split}
  \end{equation*}
  
  Next, with \eqref{eq39},  the Markov inequality, and Condition \eqref{eq23} we get for any $\epsilon >0$
  \begin{equation*}\label{th6.p5}
  \begin{split}
  &\sum_{h\in\mathcal{H}}\bar{\pi}_h \mathbf{1}\biggl\{  S\biggl(-\frac{A_h}{\sigma_h}+x\biggr)
\le  \exp\biggl[ q_\alpha^\varkappa+\log\frac{n_h}{\bar{\pi}_h}\biggr] \biggr\} \\
 & \quad \le \sum_{h\in\mathcal{H}}\bar{\pi}_h \mathbf{1}\biggl\{ r_h(q_\alpha^\varkappa,\bar{\pi}_h)-\biggl(\frac{A_{h}}{\sigma_h}-x \biggr) ^2
    > o(1) \biggr\}   
 \\
 &\quad    \le  \frac{1}{\epsilon H^*_{\bar{\pi}}}\sum_{h\in\mathcal{H}}\bar{\pi}_h \biggr[ r_h(q_\alpha^\varkappa,\bar{\pi}_h)-
     \biggl( \frac{A_{h}}{\sigma_h}-x\biggr)^2 +\epsilon H^*_{\bar{\pi}}\biggr]_+
     \\
   &\quad  \le  \frac{1}{\epsilon H^*_{\bar{\pi}}}\sum_{h\in\mathcal{H}}\bar{\pi}_h \biggr[ R_h(q_\alpha^\varkappa,H_{\bar{\pi}})-
     \biggl( \frac{A_{h}}{\sigma_h}-x\biggr)^2 +\epsilon H^*_{\bar{\pi}} \biggr]_++o(1)
     .
  \end{split}
  \end{equation*}
  To complete the proof, let us choose $x=\sqrt{\epsilon H^*_{\bar{\pi}}}$.

\subsection{Proof of Theorem \ref{th9}}

\textbf{ A lower bound.} For given $x$ and $\delta>0$  by \eqref{eq17} and \eqref{eq13},    we obtain 
\begin{equation*}
\begin{split}
&\beta^B_{\rho,\tau}(A) \ge
\mathbf{P} \biggl\{\frac{\bar{\pi}_\rho}{n_\rho}S\biggl(-\frac{A}{\sigma_\rho}+\xi_{\rho,\tau}\biggr)
+\zeta^\circ\le q_\alpha^\circ;\xi_{\rho,\tau}\ge x \biggr\}\\
& \quad \ge
\mathbf{P} \biggl\{\frac{\bar{\pi}_\rho}{n_\rho}S\biggl(-\frac{A}{\sigma_\rho}+x\biggr)
+\zeta^\circ\le q_\alpha^\circ;\xi_{\rho,\tau}\ge x \biggr\}
\\&\quad =\mathbf{P}\bigl\{\xi_{\rho,\tau}\ge x\bigr\}
\mathbf{P} \biggl\{\zeta^\circ\le q_\alpha^\circ-\frac{\bar{\pi}_\rho}{n_\rho}S\biggl(-\frac{A}{\sigma_\rho}+x\biggr)\biggr\}\\&\quad\ge 
\mathbf{P}\bigl\{\xi_{\rho,\tau}\ge x\bigr\}
\mathbf{P} \bigl\{\zeta^\circ\le (1-\delta) q_\alpha^\circ\bigr\}
\mathbf{1}\biggl\{\frac{\bar{\pi}_\rho}{n_\rho}S\biggl(-\frac{A}{\sigma_\rho}+x\biggr)\le  \delta q_\alpha^\circ \biggr\}.
\end{split}
\end{equation*}

Similarly to \eqref{eq40}, since $\bar{\pi}\rightarrow 0$, we get
\begin{equation*}
\begin{split}
&\mathbf{E}_{\bar{\pi}}\biggl[\mathbf{1}\biggl\{\frac{\bar{\pi}_h}{n_h}S\biggl(-\frac{A_h}{\sigma_h}+x\biggr)\le  \delta q_\alpha^\circ \biggr\} \biggr]
\\ &\quad
\ge 1-\frac{1}{\epsilon  H^*_{\bar{\pi}} }\mathbf{E}_{\bar{\pi}} \biggl[\biggl(
\frac{A_h}{\sigma_h}-x\biggr)^2+\epsilon H^*_{\bar{\pi}} -R_h[\log(\delta q_\alpha^\circ), H_{\bar{\pi}}]
\biggr]_+\\
&\quad
\ge 1-\frac{1}{\epsilon  H^*_{\bar{\pi}} }\mathbf{E}_{\bar{\pi}} \biggl[\biggl(
\frac{A_h}{\sigma_h}-x\biggr)^2+\epsilon H^*_{\bar{\pi}} -R_h[\log( q_\alpha^\circ), H_{\bar{\pi}}]\biggr]_+\\
&\quad \ge 1+o(1). 
\end{split}
\end{equation*}
So, since $\delta$ is arbitrary, \eqref{eq29} follows from the above inequalities.

\medskip
\noindent
\textbf{ An upper bound.} Since $S(x)$ is decreasing, we get
\begin{equation}\label{eq41}
\begin{split}
\beta^B_{\rho,\tau}(A) \le&
\mathbf{P} \biggl\{\frac{\bar{\pi}_\rho}{n_\rho}S\biggl(-\frac{A}{\sigma_\rho}+\xi_{\rho,\tau}\biggr)
+\zeta^\circ\le q_\alpha^\circ;\xi_{\rho,\tau}\ge x \biggr\}\\
 +& \mathbf{P} \biggl\{\frac{\bar{\pi}_\rho}{n_\rho}S\biggl(-\frac{A}{\sigma_\rho}+\xi_{\rho,\tau}\biggr)
+\zeta^\circ\le q_\alpha^\circ;\xi_{\rho,\tau}\le x \biggr\}\\
\le&
\mathbf{P} \bigl\{ \xi_{\rho,\tau}\ge x \bigr\}
 + \mathbf{P} \biggl\{\frac{\bar{\pi}_\rho}{n_\rho}S\biggl(-\frac{A}{\sigma_\rho}+x \biggr)
+\zeta^\circ\le q_\alpha^\circ \biggr\}.
\end{split}
\end{equation}

Next, for any given $x^\circ< q_\alpha^\circ$ we obtain with the help of the Markov inequality and \eqref{eq39}
\begin{equation*}
\begin{split}
&\mathbf{E}_{\bar{\pi}}\biggl[\mathbf{P} \biggl\{\frac{\bar{\pi}_h}{n_h}S\biggl(-\frac{A_h}{\sigma_h}+x \biggr)
+\zeta^\circ\le q_\alpha^\circ \biggr\}\biggr]\\ 
&\quad \le \mathbf{E}_{\bar{\pi}}\biggl[\mathbf{P} \biggl\{\frac{\bar{\pi}_h}{n_h}S\biggl(-\frac{A_h}{\sigma_h}+x \biggr)
+\zeta^\circ\le q_\alpha^\circ; \zeta^\circ<x^\circ \biggr\}\biggr]\\
&\qquad+\mathbf{E}_{\bar{\pi}}\biggl[\mathbf{P} \biggl\{\frac{\bar{\pi}_h}{n_h}S\biggl(-\frac{A_h}{\sigma_h}+x \biggr)
+\zeta^\circ\le q_\alpha^\circ; \zeta^\circ\ge x^\circ \biggr\}\biggr]\\
&\quad \le \mathbf{P}\bigl\{\zeta^\circ<x^\circ\} +\mathbf{E}_{\bar{\pi}}\biggl[\mathbf{1} \biggl\{\frac{\bar{\pi}_h}{n_h}S\biggl(-\frac{A_h}{\sigma_h}+x \biggr)
\le q_\alpha^\circ-x^\circ \biggr\}\biggr]\\
&\quad\le  \mathbf{P}\bigl\{\zeta^\circ<x^\circ\}+\frac{1}{\epsilon  H^*_{\bar{\pi}} }\mathbf{E}_{\bar{\pi}} \biggl[ R_h[\log( q_\alpha^\circ-x^\circ), H_{\bar{\pi}}]+\epsilon H^*_{\bar{\pi}}-\biggl(
\frac{A_h}{\sigma_h}-x\biggr)^2\biggr]_+.
\end{split}
\end{equation*}
Finally choosing $x=\sqrt{\epsilon H^*_{\bar{\pi}}}$ and combining this equation with \eqref{eq41}, we complete the proof of \eqref{eq35.x}.

\subsection{Proof of Theorem \ref{th7}}
 Let us consider 
\[
H\bigl(\bar{\pi}^{\omega,\nu},{\Pi}\bigr)=\sum_{h\in \mathcal{H}} \bar{\pi}^{\omega,\nu}_h \log\frac{1}{{\Pi}_h}.
\]
One can check easily with \eqref{eq34}--\eqref{eq34.2}  and \eqref{prior} that as $\omega\rightarrow\infty$
\begin{equation*}
\begin{split}
&H\bigl(\bar{\pi}^{\omega,\nu},{\Pi}\bigr)=\log\frac{1}{\varepsilon}+\frac{1}{\omega}\sum_{k=1}^\infty 
\nu\biggl(\frac{k}{\omega}\biggr) L^{m,\varepsilon}(k)+o(1) =
\log\frac{1}{\varepsilon}+o(1)\\ &+\int_0^\infty \nu(x)\biggl[\log(x\omega+1)+\sum_{s=0}^{m-1}\log[\psi_s(x\omega)]+(1+\varepsilon)\log[\psi_m(x\omega)]\biggr]\, dx.
\end{split}
\end{equation*} 
It is also clear in view of \eqref{eq31} that 
\[
\int_0^\infty \nu(x)\log(x\omega+1)\, dx = \log(\omega)+O(1) 
\]
and for any integer $s\ge 0$
\[
\int_0^\infty \nu(x)\log[\psi_s(x\omega)]\, dx = \log[\psi_s(\omega)]+O(1). 
\]
Therefore as $\omega\rightarrow\infty$
\begin{equation}\label{eq42}
H\bigl(\bar{\pi}^{\omega,\nu},{\Pi}\bigr)= \log(\omega)+(1+o(1))\log^*(\omega). 
\end{equation}

Next, with  similar arguments we obtain
\begin{equation}\label{eq43}
\begin{split}
&\sum_{h\in \mathcal{H}}\bar{\pi}^{\omega,\nu}_h \biggl[H(\bar{\pi}^{\omega,\nu},{\Pi})
-\log\frac{1}{{\Pi}_h}\biggr]_+\\ &\qquad
=\int_{0}^\infty \nu(x)\Bigl[H(\bar{\pi}^{\omega,\nu},{\Pi})-\log(x\omega+1)\\ &\qquad \quad -(1+o(1))\log[\log(x\omega+1)+1] \Bigr]_+dx\\ &\qquad=
\int_{0}^\infty \nu(x)\Bigl[\log(\omega)+(1+o(1))\log^*(\omega)]-\log(x\omega+1)\\ &\qquad \quad -(1+o(1))\log[\log(x\omega+1)+1] \Bigr]_+dx\\ &\qquad
=
\int_{0}^\infty \nu(x)\bigl[-\log(x+1)+o(1)\log^*(\omega)\bigr]_+dx =o(1)\log^*(\omega)].
\end{split}
\end{equation}

 In view of \eqref{eq42} and \eqref{eq43} the rest of the proof is similar to the one of Theorem \ref{th6} and therefore omitted.

\end{document}